\documentclass[11pt]{amsart} 
\usepackage{amsmath,latexsym,graphicx}
\usepackage{amssymb,amsthm,amsfonts}

\newtheorem{theorem}{Theorem}[section]
\newtheorem{lemma}[theorem]{Lemma}
\newtheorem{remark}[theorem]{Remark}
\newtheorem{proposition}[theorem]{Proposition}

\def\O{\Omega}

\def\cal{\mathcal}
\def\H{H^1_0(\Omega)}

\begin{document}

\title[Multiple positive solutions for a  SPS system ]
{Multiple positive solutions for a Schr{\"o}dinger-Poisson-Slater system }

\author{Gaetano Siciliano}
\address{Gaetano Siciliano, Dpto. An\'alisis Matem\'atico, Universidad de Granada, 18071 Granada, Spain}
\thanks{The author is supported by J. Andaluc\'{\i}a - FQM 116 and by 
M.I.U.R. - P.R.I.N. ``Metodi variazionali e topologici
nello studio di fenomeni nonlineari''.}
\email{sicilia@ugr.es}


\begin{abstract}
In this paper we investigate the existence of positive solutions to
the following Schr{\"o}dinger-Poisson-Slater system 
\begin{equation*}\left\{
\begin{array}
[c]{ll}
- \Delta u+ u + \lambda\phi u=|u|^{p-2}u & \text{ in } \Omega\\
-\Delta\phi= u^{2} & \text{ in } \Omega\\
u=\phi=0 & \text{ on } \partial\Omega
\end{array}
\right.
\end{equation*}
where $\Omega$ is a bounded domain in $\mathbf{R}^{3},\lambda$ is a fixed positive parameter and $p<2^{*}=\frac{2N}{N-2}$.
We prove that if $p$ is ``near'' the critical Sobolev exponent $2^*$, then the number of positive solutions
is greater then 
the Ljusternik-Schnirelmann category of $\Omega$.
\end{abstract}
\maketitle

\section{Introduction}

In  \cite{BC,BC2} Benci and Cerami proved a result on the number
of positive solutions of the following problem
\begin{equation}
\label{model}\left\{
\begin{array}
[c]{ll}
-\Delta u +u=|u|^{p-2}u & \quad\text{ in } \Omega\\
u=0 & \quad\text{ on } \partial\Omega
\end{array}
\right.
\end{equation}
where $\Omega\subset\mathbf{R}^{N}$ is a smooth and bounded domain, $N\ge3$
and $p<2^{*}=\frac{2N}{N-2}$, the critical Sobolev exponent for the embedding
of $\H$ in $L^p(\Omega)$. In particular
they ask how the number of positive solutions depends on the topology of
$\Omega$. The core of their results is that if $\Omega$ is ``topologically
rich''
then there are many solutions as soon as the nonlinearity acts strongly on the
equation. For problem \eqref{model} this happens when $p$ is near $2^{*}$;
indeed they prove the following result

\begin{theorem}
\label{ThA} There exists a $\bar{p}\in(2,2^{*})$ such that for every
$p\in[\bar p, 2^{*})$ problem \eqref{model} has (at least) $\emph{cat}_{\bar\Omega}\,(\bar\Omega)+1$
positive solutions.
\end{theorem}

Hereafter $cat$ is the Ljusternik-Schnirelmann category (see e.g. \cite{J}).

They prove Theorem \ref{ThA} by variational methods 
looking for the solutions as critical points 
of an energy functional restricted to a suitable manifold on which it is bounded from below.
Then, since the Palais-Smale
condition  (see below for the definition)  is satisfied 
the main effort is to found a sublevel of the functional with a non-zero category,
let us say $k$; in these conditions the Ljusternik-Schnirelmann theory would give the existence of at least $k$ critical points.
By  introducing the barycenter map, they are able to find sublevels
with category greater then the category of $\Omega$ and so
the existence of at least $\mbox{cat}_{\bar\Omega}\,(\bar\Omega)$
critical points is ensured. Actually this is done in \cite{BC} while the existence of another solution 
is proved in \cite{BC2}.

 Another approach with the Morse theory has been used in \cite{BC3} for more general nonlinearity than $|u|^{p-2}u.$

We need to recall that problems like \eqref{model}, in bounded or exterior domain, even with the critical exponent 
and with a control parameter $\varepsilon>0$ have  been object of wide investigation. Also the concentration (blow-up)
of solutions in specific points of the domain $\Omega$ when the parameter tends to zero is studied: we limit  ourselves 
to citing \cite{CP, DF, P, R2, W} and the references therein.

\medskip

The aim of this paper is to prove an analogous result of Theorem \ref{ThA} for
the Schr{\"o}dinger-Poisson-Slater system:
\begin{equation}
\label{SPS}\left\{
\begin{array}
[c]{ll}
- \Delta u+ \omega u + \lambda\phi u=|u|^{p-2}u & \text{ in } \Omega,\\
-\Delta\phi= u^{2} & \text{ in } \Omega,\\
u=\phi=0 & \text{ on } \partial\Omega,
\end{array}
\right.
\end{equation}where $\Omega$ is a (smooth and) bounded domain in $\mathbf{R}^{3}$,
$p\in(2,2^{*}), \omega>0$ and $\lambda$ is a positive fixed parameter. It is assumed 
$\mbox{cat}_{\bar\Omega}\,(\bar\Omega)>1$.

This system appears studying 
the nonlinear Schr{\"o}dinger equation
$$i\hbar\frac{\partial \psi}{\partial t}=-\frac{\hbar^2}{2m}\Delta_x \psi+|\psi|^{p-2}\psi$$ 
which describes   quantum (non relativistic) particles  interacting  with the electromagnetic field generated
by the motion. Here $\psi=\psi(x,t)$ is a complex valued function and $\hbar, m>0$ are interpreted respectively as
the normalized Plank constant and the mass of the particle. However, since they have no role in our analysis,
we set $\hbar=1$ and $m=1/2$. A model for the interaction between matter and electromagnetic 
field is provided by the abelian gauge theories
but can also be derived by the Slater  approach to the  Hartree-Fock model. Without entering in details
(the reader interested is refereed e.g. to \cite{Bleecker,slater}),
if $\phi(x,t)$ and $\mathbf A(x,t)$ denote the gauge potentials
of the e.m. field, the search of stationary solutions, namely solutions
$\psi$ of the form 
$$\psi(x,t)=u(x)e^{i\omega t}\ \ \ u(x)\in \mathbf R,\,\omega>0\,,$$
in the purely electrostatic case
$$\phi=\phi(x) \ \ \text{and}\ \ \mathbf A=\mathbf 0\,,$$
leads exactly to the system we want to study. The boundary conditions $u=\phi=0$ on
$ \partial\Omega $ mean that the particle is constraint to live in $ \Omega $.
In the following, referring to \eqref{SPS} we will assume for simplicity $\omega=1$.

Problem \eqref{SPS} contains two kinds of nonlinearities: the first one is $\phi u$ and concerns 
the interaction with the electric field.
This nonlinear term is nonlocal since the electrostatic potential $\phi$ depends also on 
the wave function to which is related  by the Poisson equation $-\Delta\phi= |\psi|^2=u^{2}$.
The second nonlinearity is $|u|^{p-2}u$ . This one contains the Slater 
correction term $C_S\, |u|^{2/3}u$, where $C_S$ is the
Slater constant and depends on the particles considered (for more details see
\cite{BoLoSo, slater}). Physically speaking, the local nonlinearity $|u|^{p-2}u$ represents
the interaction among many particles and is in competition with the  intrinsic nonlinearity
of the system $\phi u$.

Motivated by some perturbation results
(see e.g. \cite{wei,m3as} in which the case with $\Omega=\mathbf{R}^3$  and $\lambda\rightarrow0^+$ is considered),
we have introduced the parameter $\lambda>0$ which takes a role also
in a bounded domain, at least for small values of $p$. 

Because of its importance in many different physical framework, the Schr{\"o}dinger-Poisson-Slater system 
(sometimes called Schr{\"o}dinger-Maxwell system)
has been extensively studied in the past
years: besides the results on bounded domains (see e.g. \cite{BF, PS1,
PS2, RS}), there are also many papers on $\mathbf{R}^{3}$ which treat different
aspects of the SPS system, even with an additional external and fixed potential $V(x)$.
In particular ground states, radially and non-radially solutions
or semiclassical limit and concentration of solutions are studied, see e.g.
 \cite{aruiz, AzzPo, Coc, DApMu, dAv, Kik, IV, jfa, WZh}.

\medskip

We approach problem \eqref{SPS}  by variational methods: the weak
solutions  are characterized as critical points of a 
$C^{1}$ functional $I= I(u)$ defined on the Sobolev space $H^{1}_{0}(\Omega)$ or
a suitable submanifold (see below). 
A fundamental tool to apply variational techniques is the
so-called \emph{Palais-Smale condition} (PS for brevity): every sequence
$\{u_{n}\}$ such that
\begin{equation}\label{PS}
\{I(u_{n})\} \ \text{ is bounded \ \  and }\ \ I^{\prime}(u_{n})\rightarrow0 \ \ \text{ in } H^{-1}(\Omega),
\end{equation} 
admits a converging subsequence. Sequences which satisfy
\eqref{PS} are called \emph{Palais-Smale sequences.}

Now, it is known that when
$p\in(4,2^*)$ the PS condition holds (see e.g. \cite{PS2}), hence we have hope
to apply classical theorems of LS theory in the same spirit of \cite{BC} and \cite{BC2}, to find critical points of $I$; 
indeed we get the following result
\begin{theorem}\label{Main}
There exists a $\bar{p}\in(4,2^*)$ such that for every $p\in[\bar{p},2^*)$ problem
\eqref{SPS} has at least $\emph{cat}_{\bar{\Omega}}\,(\bar{\Omega})+1$ 
positive solutions.
\end{theorem}
It is understood that $\bar{p}$ does not depend on the ``strength" of the interaction $\lambda.$
We remark that the weak solutions  found by means of the
variational method are indeed classical solutions, by standard regularity
results.

To prove the theorem we use the general ideas of Benci and Cerami adapting 
their arguments to our problem which contains also the coupling term $\phi u$.

\medskip

The paper is organized
as follow: in the next Section we fix the notations and recall some
useful facts. Sections 3 and 4 are devoted to the functional setting and
to introduce the ingredients which allow us to use the abstract theory
of Ljusternik-Schnirelmann. Finally the proof of Theorem \ref{Main}
is completed in Section 5.

\subsubsection*{Acknowledgment}
The author wishes to thank David Ruiz for helpful discussions on the matter and prof. G. Cerami
for bringing to his attention paper \cite{BC2}.

\section{Some notations and preliminaries}

Without loss of generality we assume in all the paper $0\in\Omega$. We denote by
$|\,.\,|_{L^{p}(A)}$ the $L^{p}-$norm of a function defined on the domain $A$.
If the domain is specified (usually $\Omega$) or if there is no confusion, we use the notation
$|\,.\,|_{p}$. Moreover let $H^{1}_{0}(\Omega)$ be the usual Sobolev space
with (squared) norm
$$
\|u\|^{2}=|\nabla u|_{2}^{2}+|u|_{2}^{2}%
$$
and dual $H^{-1}(\Omega)$.

We use $B_{r}(y)$ for the closed ball of radius $r>0$ centered in $y$. If $y=0$
we simply write $B_r.$

The letter $c$ will be used indiscriminately
to denote a suitable positive constant whose value may change from line to line 
and we will use $o(1)$ for a quantity which goes to zero.

Finally, in view of our Theorem \ref{Main}, from now on we assume $p>4$.
Other notations will be introduced in Section \ref{bary}.

\medskip

First of all, let $\phi_{u}\in H^{1}_{0}(\Omega)$ be the unique (and positive) solution of
$-\Delta\phi=u^{2}$ and $\phi=0$ on $\partial\Omega$ and let us recall the
following properties that will be repeatedly used (for a proof see e.g. \cite{jfa}):

\begin{itemize}
\item for any $\alpha,\beta\ge0,t>0$ let $u_{t}(\cdot)=t^{\alpha}u(t^{\beta
}(\cdot))$. Then
\[
\phi_{u_{t}}(\cdot)=t^{2(\alpha-\beta)}\phi_{u} (t^{\beta}(\cdot))\,;
\]

\item $u_{n}\rightharpoonup u$ in $H^{1}_{0}(\Omega)\Longrightarrow
\int_{\Omega} \phi_{u_{n}} u_{n}^{2}dx \rightarrow\int_{\O}\phi_{u} u^{2}\,dx$\,; \medskip

\item $|\nabla\phi_{u}|_{2}\le c|\nabla u|^2_{2}$ for some constant
$c>0$\,;\medskip

\item $\int_{\Omega}|\nabla\phi_{u}|^{2}dx=\int_{\O}\phi_{u} u^{2}dx$\,.
\end{itemize}

\medskip

The functional associated to \eqref{SPS} is
\begin{equation}
\label{I}I_{p}(u)=\frac{1}{2}\int_{\O}(|\nabla u|^{2}+u^{2})\,dx+\frac{\lambda}{4}%
\int_{\O}\phi_{u} u^{2}dx-\frac{1}{p}\int_{\O}|u|^{p} dx
\end{equation}
and its critical points are the solutions of the system 
(see e.g. \cite{BF}). However the functional is unbounded from above and from below on
$H^{1}_{0}(\Omega)$. The idea is to restrict the functional to  a suitable manifold on 
which this unboundedness is removed.

\medskip

In \cite{BC2} the authors deal with $E(u)=\frac{1}{2}\|u\|^2-\frac{1}{p}|u|_p^p$
and
to overcome the unboundedness they introduce the constraint
$$V_{p}=\left\{u \in H^{1}_{0}(\Omega) : |u|_{p}=1 \right\}.$$
On $V_p$ the functional $E$ is bounded from below (achieves its minimum), satisfies
the PS condition and the
classical LS theory applies. This gives constraint critical points and
Lagrange multipliers appear in the right hand side of the equation in (\ref{model}).
Finally, ``stretching'' the multipliers one
gets solutions of (\ref{model}).

\medskip

In our case the constraint $V_{p}$ is not a good choice although $I_{p}$ would
have a minimum on $V_{p}.$ This is due to a different degree of homogeneity of
the added term $\lambda\phi_{u} u$; indeed  it is easy to see that
there is no way to eliminate the Lagrange multiplier once it appears.
We study the functional \eqref{I} on a \emph{natural constraint} and
in this case the Nehari manifold works well.

\section{The Nehari manifold}

In this section we recall some known facts 
about the Nehari manifold  that will be used throughout the paper. %

The Nehari manifold associated to \eqref{I} is defined
by
\begin{equation*}
\label{Nehari}{\mathcal{N}}_{p}=\left\{ u\in H^{1}_{0}(\Omega)\setminus\{0\}:
G_{p}(u)=0\right\}
\end{equation*}
where
$$
G_{p}(u):=I_{p}^{\prime}(u)[u]=\|u\|^{2}+\lambda\int_{\Omega}\phi_{u}
u^{2}dx-|u|_{p}^{p}\,.
$$
On $\mathcal{N}_{p}$ the functional \eqref{I} has the form
\begin{equation}
\label{Ivinc}I_{p}(u)=\frac{p-2}{2p}\|u\|^{2}+\lambda\frac{p-4}{4p}
\int_{\Omega}\phi_{u} u^{2} dx\,.
\end{equation}
Sometimes we will refer to \eqref{Ivinc} as the constraint functional, also
denoted with $I_{p}|_{\mathcal{N}_{p}}$. 

In the next Lemma we recall the basic properties of the Nehari manifold. 
\begin{lemma}
\label{lemmanehari}
We have
\begin{itemize}
\item[1.] $\mathcal{N}_{p}$ is a $C^{1}$ manifold\,,

\item[2.] there exists $c>0$ such that for every $u\in \mathcal{N}_{p}:c\le\|u\|\,, $

\item[3.] for every $u\neq0$ there exists a unique $t>0$ such that
$tu\in\mathcal{N}_{p}$\,,

\item[4.] the following equalities are true
\[
m_{p}=\inf_{u\neq0}\max_{t>0}I_{p}(tu)=\inf_{g\in\Gamma_{p}} \, \max
_{t\in[0,1]} I_{p}(g(t))
\]
where
$$
\Gamma_{p}=\left\{g\in C([0,1];H^{1}_{0}(\Omega)) : g(0)=0, I_{p}(g(1))\le0,
g(1)\neq0\right\}.
$$
\end{itemize}
\end{lemma}

Then recalling that $p>4$, we have 
$$m_{p}:=\inf_{u\in\mathcal{N}_{p}}I_{p}(u)>0\,.$$
Moreover the manifold ${\mathcal{N}}_{p}$ is a natural constraint for $I_{ p}$ (given by \eqref{I}) in the
sense that any $u\in\mathcal{N}_{p}$ critical point of $I_{p}|_{\mathcal{N}%
_{p}}$ is also a critical point for the free functional $I_{p}$ (for a proof of
these facts, see e.g. Section 6.4
in \cite{AM}). Hence the
(constraint) critical points we find are solutions of our problem since no
Lagrange multipliers appear.

The Nehari manifold well-behaves with respect to the PS sequences:

\begin{lemma}\label{lemmaPS}
Let $\{u_{n}\}\subset\mathcal{N}_{p}$ be a PS sequence
for $I_{p}|_{\mathcal{N}_{p}}$. Then it is a PS
sequence for the free functional $I_{p}$ on $H^{1}_{0}(\Omega)$.
\end{lemma}
\begin{proof}
By definition, $\{u_{n}\}\subset\mathcal{N}_{p}$, $I_{p}|_{\mathcal{N}_{p}}(u_{n})$ 
is bounded and there exist Lagrange multipliers $\{\mu_{n}%
\}\subset\mathbf{R}$ such that $(I_{p}|_{\mathcal{N}_{ p}})^{\prime
}(u_{n})=I^{\prime}_{ p}(u_{n})-\mu_{n} G^{\prime}_{ p}(u_{n})\rightarrow0$ in
$H^{-1}(\Omega)$. Then recalling the definition of $\mathcal{N}_{ p}$ we have
$$(I_{p}|_{\mathcal{N}_{p}})^{\prime}(u_{n})[u_n]=\mu_{n} G^{\prime}_{ p}(u_{n})[u_{n}]\rightarrow0.$$
Since $G^{\prime}_{ p}(u_{n})[u_{n}]\neq0$ it follows that the sequence of
multipliers vanishes and
\begin{equation*}
I^{\prime}_{ p}(u_{n})=(I_{ p}|_{\mathcal{N}_{ p}})^{\prime}(u_{n})+\mu_{n}
G^{\prime}_{ p}(u_{n})\rightarrow0.
\end{equation*}
\end{proof}

As we have already anticipated, for $p\in(4,2^{*})$ it is known that the free functional $I_{ p}$ given by
\eqref{I} satisfies the PS condition on $H^{1}_{0}(\Omega)$ (see e.g.
\cite{PS2}). 
The fact that the PS condition follows also for the
functional restricted to $\cal N_p$ is standard. 

In the following we will deal always with the restricted functional on the Nehari manifold; this
will be denoted simply with $I_{p}.$

As a consequence of the PS condition we deduce that
\begin{equation*}\label{mp}
\forall\,p\in(4,2^{*})\,:\ \ m_{ p}=\min_{\mathcal{N}_{ p}} I_{p}=I_{ p}(u_{ p})\,,
\end{equation*}
i.e. $m_{ p}$ is achieved on a function, hereafter denoted with $u_{ p}$, in $\cal N_p$. 
Since $u_{ p}$ minimizes the energy
$I_p$, it will be called a \emph{ground state}.

Observe that the sequence of minimizers $\{u_{ p}\}_{p\in(4,2^{*})}$ is bounded away from
zero; indeed, since $u_{ p}\in\mathcal{N}_{ p}\,$,
\begin{equation}
\label{dis}\|u_{ p}\|^{2}\le|u_{ p}|_{p}^{p}\le C\|u_{ p}\|^{p}
\end{equation}
where $C$ is a positive constant which can be made independent of $p$.
Hence
\begin{equation*}
\exists\,c>0 \ \ \mbox{ s.t. } \forall\, p\in(4,2^{*})\, :\ \ 0< c\le\|u_{
p}\|.
\end{equation*}

\begin{remark}
\label{rem} Turning back to \eqref{dis}, we have that $\{|u_{ p}|_{p}%
\}_{p\in(4,2^{*})}$ is bounded away from zero. Moreover, denoting with $|\Omega|$ the
Lebesgue measure of $\Omega$, by the H\"{o}lder
inequality,
\[
|u_{ p}|_{p}\le|\Omega|^{\frac{2^{*}-p}{2^{*}p}}|u_{ p}|_{2^{*}}%
\]
and so also $\{|u_{ p}|_{2^{*}}\}_{p\in(4,2^{*})}$ is bounded away from zero.
\end{remark}

Clearly, all we have stated until now is true also in the case $\lambda=0$. Moreover
also the case $p=2^*$ is covered for those results which do not require compactness
(in particular Lemma \ref{lemmanehari} and \ref{lemmaPS}).

\subsection{The limit cases}
We consider in this subsection two limit cases related to \eqref{SPS}. Our intent is to evaluate
the limit of the sequence $\{m_p\}_{p\in(4,2^*)}$
when $p\rightarrow2^*$.

\medskip

The first case is the critical problem. Let us introduce the functional
\[
I_{*}(u)=\frac{1}{2}\|u\|^{2}-\frac{1}{2^{*}}|u|_{2^{*}}^{2^{*}}\,
\]
whose critical points are the solutions of
\begin{equation}
\left\{
\begin{array}
[c]{ll}%
\label{star} -\Delta u +u=|u|^{2^{*}-2}u & \quad\text{ in } \Omega\\
u=0 & \quad\text{ on } \partial\Omega.
\end{array}
\right.
\end{equation}
It is known that the lack of compactness of the embedding of $H^{1}_{0}(\Omega)$
in $L^{2^{*}}(\Omega)$ implies that 
$I_{*}$ does not satisfies the PS condition at every level.
This is due to the invariance with respect to the conformal scaling
$$u(\cdot)\longmapsto u_{R}(\cdot):=R^{1/2}u(R(\cdot)) \ \ \ \ (R>1)$$
which leaves invariant the $L^2-$norm of the gradient an the $L^{2^*}-$norm, i.e.
$|\nabla u_{R}|_{2}^{2}=|\nabla u|_{2}^{2}$ and $|u_{R}|_{2^{*}}^{2^{*}}=|u|_{2^{*}}^{2^{*}}\,.$

As a consequence, if
$$
\mathcal{N}_{*}=\{u\in H^{1}_{0}(\Omega): G_{*}(u)=0\}\,, \ \ \ G_{*}
(u)=\|u\|^{2}-|u|_{2^{*}}^{2^{*}}%
$$
is the Nehari manifold associated, it can be proved that
$$
m_{*}:=\inf_{\mathcal{N}_{*}} I_{*} \ \ \mbox{ is not achieved. }
$$
The following lemma is known but for the sake of completeness we give the proof.
\begin{lemma}
There holds
$$
m_{*}=\frac{1}{3}S^{3/2}$$
where $S=\inf_{u\in H^{1}_{0}
(\Omega), u\neq0} \frac{\|u\|^{2}}{|u|_{2^{*}}^{2}}$ is the best Sobolev constant.
\end{lemma}
\begin{proof}
This is indeed an easy computation. First observe that for $A,B>0$ it results
\[
\max_{t>0} \left\{ \frac{t^{2}}{2}A-\frac{t^{2^{*}}}{2^{*}}B\right\} =\frac
{1}{3}\left( \frac{A}{B^{1/3}}\right) ^{3/2}.
\]
Then
\[
m_{*}=\inf_{u\neq0}\max_{t>0} I_{*}(tu)=\frac{1}{3}\left( \inf_{u\neq0}
\frac{\|u\|^{2}}{|u|_{2^{*}}^{2}}\right) ^{3/2} =\frac{1}{3}S^{3/2}.
\]
\end{proof}

The value $m_*$ turns out to be an upper bound for the sequence of ground states levels $\{m_p\}_{p\in(4,2^*)}$.
Before to prove this, let us observe that, as easy computations show:
\begin{enumerate}
\item[1)] $ |u_{R}|_{p}^{p}=R^{\frac{p-2^{*}}{2}}|u|_{p}^{p}\,,$ \medskip
\item[2)] $\int_{\Omega}\phi_{u_{R}} u_{R}^{2}\,dx=R^{-3}\int_{\Omega}\phi_{u}u^{2}\,dx\,.$\medskip
\end{enumerate}

\begin{lemma}\label{limitate}
We have
$$\limsup_{p\rightarrow2^{*}} m_{p}\le m_{*}.$$
\end{lemma}
\begin{proof}
Fix $\varepsilon>0.$ By definition of $m_{*}$ there exists $u\in
\mathcal{N}_{*}$ such that
\begin{equation}
\label{e/2}I_{*}(u)=\frac{1}{2}\|u\|^{2}-\frac{1}{2^{*}}|u|_{2^{*}}^{2^{*}
}=\frac{1}{3}\|u\|^{2}<m_{*}+\frac{\varepsilon}{2}.
\end{equation}
For $R>1$ (to be specified later), we have
\[
I_{*}(u_{R})=\frac{1}{2}|\nabla u|_{2}^{2}+\frac{1}{2R^{2}}|u|_{2}^{2}
-\frac{1}{2^{*}}|u|_{2^{*}}^{2^{*}}<m_{*}+\frac{\varepsilon}{2}\,.
\]
Now consider, for any $p\in(4,2^{*})$, the unique positive value $t_{p}$ such
that $t_{p} u_{R}\in\mathcal{N}_{p}\,.$ 
By definition, $t_{p}$ satisfies
\begin{equation}
\label{rel}\|t_{p} u_{R}\|^{2}+\lambda t_{p}^{4} \int_{\Omega}\phi_{u_{R}}
u_{R}^{2}\,dx=|t_{p} u_{R}|_{p}^{p}
\end{equation}
from which we deduce:

\begin{itemize}
\item $\{t_{p}\}_{p\in(4,2^{*})}$ is bounded away from zero.
\end{itemize}
Indeed by \eqref{rel} and the embedding of $L^{p}$ in $H^{1}_{0}$ we get
$\|t_{p} u_{R}\|^{2}\le C\|t_{p} u_{R}\|^{p}$ so $\|t_{p} u_{R}\|^{2}\ge c$
and finally $t_{p}^{2}\ge\frac{c}{\|u_{R}\|^{2}}\ge\frac{c}{\|u\|^{2}}>0.$

\begin{itemize}
\item $\{t_{p}\}_{p\in(4,2^{*})}$ is bounded above.
\end{itemize}
Indeed $$\frac{\|u_{R}\|^{2}}{t_{p}^{2}}+\lambda\int_{\Omega}\phi_{u_{R}}
u_{R}^{2} dx =t_{p}^{p-4}|u_{R}|_{p}^{p}$$ and, by the continuity of the map
$p\mapsto|u_{R}|_{p}$\,, it is readily seen that if $t_{p}$ tends to
$+\infty$ we get a contradiction.

\smallskip

So we may assume that $\lim_{p\rightarrow2^{*}} t_{p}=t_{*}$ and passing to
the limit in \eqref{rel} we get
\begin{align*}
t_{*}^{2}|\nabla u|_{2}^{2}+\frac{t_{*}^{2}}{R^{2}}|u|_{2}^{2}+\lambda
\frac{t_{*}^{4}}{R^{3}}\int_{\Omega}\phi_{u} u^{2}dx & =t_{*}^{2^{*}%
}|u|_{2^{*}}^{2^{*}}\\
& =t_{*}^{2^{*}}\left( |\nabla u|_{2}^{2}+|u|_{2}^{2}\right)
\end{align*}
or equivalently,
\[
(t_{*}^{2^{*}}-t_{*}^{2})|\nabla u|_{2}^{2}=\frac{t_{*}^{2}}{R^{2}}|u|_{2}%
^{2}+\lambda\frac{t_{*}^{4}}{R^{3}}\int_{\Omega}\phi_{u} u^{2}dx-t_{*}^{2^{*}%
}|u|_{2}^{2}\,.
\]
Now if $R$ is chosen sufficiently large, the r.h.s. above is negative and 
we deduce
\begin{equation}
\label{t}t_{*}<1.
\end{equation}
Furthermore
\begin{eqnarray*}
I_{p}(t_{p} u_{R})&=&\frac{p-2}{2p}\|t_{p} u_{R}\|^{2}
+\lambda\frac{p-4}{4p}t_{p}^{4}\int_{\Omega}\phi_{u_{R}} u_{R}^{2}\,dx\\
&=&\frac{p-2}{2p}t^{2}_{p}|\nabla u|^{2}_{2}+\frac{p-2}
{2p}\frac{t_{p}^{2}}{R^{2}}|u|_{2}^{2}+\lambda\frac{p-4}{4p}\frac{t_{p}^{4}
}{R^{3}}\int_{\Omega}\phi_{u} u^{2}dx 
\end{eqnarray*}
and passing to the limit for $p\rightarrow2^*$, taking advantage of \eqref{t},
\begin{align*}
\lim_{p\rightarrow2^{*}} I_{p}(t_{p} u_{R}) & =\frac{1}{3}t_{*}^{2}|\nabla
u|_{2}^{2}+\frac{1}{3}\frac{t_{*}^{2}}{R^{2}}|u|_{2}^{2}+ \frac{\lambda
t_{*}^{4}}{12 R^{3}}\int_{\Omega}\phi_{u} u^{2}\,dx\\
& <\frac{1}{3}\| u\|^{2}+\frac{\lambda}{12 R^{3}}\int_{\Omega}\phi_{u}
u^{2}\,dx\,.
\end{align*}
Lastly, if $R$ is such that $\frac{\lambda}{12
R^{3}}\int_{\Omega}\phi_{u} u^{2}\,dx<\varepsilon/2$
we get, using \eqref{e/2}
\[
\limsup_{p\rightarrow2^{*}} m_{p}\le\lim_{p\rightarrow2^{*}} I_{p} (t_{p}
u_{R})<\frac{1}{3}\|u\|^{2}+\frac{\varepsilon}{2} <m_{*}+\varepsilon
\]
which concludes the proof since $\varepsilon$ is arbitrary.
\end{proof}

Note that by \eqref{Ivinc}, the boundedness of 
$\{m_{p}\}_{p\in(4,2^{*})}$ implies the
boundedness of the ground state solutions, namely
\begin{equation}
\label{bound}\exists\, c>0 \ \ \mbox{ such that }\ \ \forall\, p\in(4,2^{*}):
\|u_{p}\|\le c.
\end{equation}

We need now a technical lemma.
\begin{lemma}\label{limsuptp}
Let $p\in(4,2^{*})$ and $t_{p}>0$ the unique value such that
$t_{p} u_{p}\in\mathcal{N}_{*}$.
Then
$$\limsup_{p\rightarrow2^{*}} t_{p}\le1.$$
\end{lemma}
\begin{proof}
By definition of $\mathcal{N}_{*}, t_{p}$ satisfies
\[
t_{p}^{2^{*}}|{u}_{p}|_{2^{*}}^{2^{*}}=t_{p}^{2}\|{u}_{p}\|^{2}%
\]
and using that ${u}_{p}\in\mathcal{N}_{p}$ and the H\"{o}lder inequality we
get
\begin{equation}
\label{tsup}t_{p}^{2^{*}-2}=\frac{|{u}_{p}|_{p}^{p}-\lambda\int_{\Omega}
\phi_{u_{p}} u_{p}^{2}\,dx}{|{u}_{p}|_{2^{*}}^{2^{*}}}\le\frac{|{u}_{p}
|_{p}^{p}}{|{u}_{p}|_{2^{*}}^{2^{*}}}\le\frac{|\Omega|^{\frac{2^{*}-p}{2^{*}}
}}{|{u}_{p}|_{2^{*}}^{2^{*}-p}}.
\end{equation}
By the embedding $L^{2^{*}
}(\Omega)\hookrightarrow H^{1}_{0}(\Omega)$ and \eqref{bound} we 
deduce that the sequence $\{|{u}_{p}|_{2^{*}}\}_{p\in(4,2^{*})}$
is bounded. 
Moreover recalling Remark \ref{rem}
we have that it is 
also bounded away from zero. So the conclusion follows by \eqref{tsup} since
$\lim_{p\rightarrow2^*}\frac{|\Omega|^{\frac{2^{*}-p}{2^{*}}}}{|{u}_{p}|_{2^{*}}^{2^{*}-p}}=1\,.$
\end{proof}
\begin{remark}\label{rembound}
Again note that Proposition \ref{limitate}, (\ref{bound}) and Lemma \ref{limsuptp}
hold also for problem \eqref{SPS} with $\lambda=0$.
\end{remark}

\smallskip

The other limit case we consider 
is that related to problem \eqref{model}, namely setting $\lambda=0$ in \eqref{SPS}.

For any $p\in(4,2^{*})$ let $\tilde I_{p}(u)=\frac{1}{2}\|u\|^{2}-\frac{1}%
{p}|u|_{p}^{p}$ be the functional on $H^{1}_{0}(\Omega)$ whose critical points
solve
\begin{equation*}
\left\{
\begin{array}[c]{ll}
- \Delta u+ u =|u|^{p-2}u & \text{ in } \Omega,\\
u=0 & \text{ on } \partial\Omega.
\end{array}
\right.
\end{equation*}
As usual, we can define $\tilde{\mathcal{N}}_{p}=\{u\in H^{1}_{0}%
(\Omega)\setminus\{0\}: \|u\|^{2}=|u|_{p}^{p}\}$ on which the functional is
$\tilde I_{p}(u)=\frac{p-2}{2p}\|u\|^{2}$ and we denote with
$$
\tilde m_{p}:=\min_{\tilde{ \mathcal{N}}_{p}} \tilde I_{p}=\tilde I_{p}%
(\tilde{u}_{p}).
$$
By Remark \ref{rembound} we have
\begin{equation}
\label{tildebound}\{\|\tilde{u}_{p}\|\}_{p\in(4,2^{*})} \ \ \mbox{is bounded.}
\end{equation}
Moreover, if
$t_{p}>0$ is such that $t_{p} u_{p}\in\tilde{\mathcal{N}_{p}},$ by \eqref{dis}
we get
$
t_{p}^{p-2}=\frac{\|u_{p}\|^{2}}{|u_{p}|_{p}^{p}}\le1$
and so
\[
\tilde{m}_{p}\le\tilde I_{p}(t_{p} u_{p})=\frac{p-2}{2p}t_{p}^{2}\|u_{p}
\|^{2}\le\frac{p-2}{2p}\|u_{p}\|^{2}< I_{p}(u_{p}).
\]
This means
\begin{equation}
\label{mm}\tilde{m}_{p}< m_{p}\,.
\end{equation}

\medskip

Now we are ready to compute the limit of $m_{p}$ when $p$ tends to $2^{*}$.

\begin{proposition}
\label{limitmp} For any bounded domain we have
\[
\lim_{p\rightarrow2^{*}} m_{p}=m_{*}.
\]
\end{proposition}

\begin{proof}
By \eqref{mm} and Lemma \ref{limitate} it is sufficient to prove that
$$m_{*}\le\liminf_{p\rightarrow2^{*}}\tilde{m}_{p}\,.$$
Let $t_{p}>0$ the unique value such that $t_{p} \tilde{u}_{p}\in
\mathcal{N}_{*}$. Applying Lemma \ref{limsuptp} (with $\lambda=0$) we know
$$\limsup_{p\rightarrow2^{*}}t_{p}\le1.$$
Finally, using \eqref{tildebound} we derive
\begin{align*}
m_{*} &  \le I_{*}(t_{p}\tilde{u}_{p})=\left( \frac{1}{2}-\frac{1}{2^{*}
}\right) t_{p}^{2}\|\tilde{u}_{p}\|^{2}\\
&  = \tilde{I}_{p}(\tilde{u}_{p}) t_{p}^{2}+\left( \frac{1}{p}-\frac{1}{2^{*}
}\right) \|\tilde{u}_{p}\|^{2}t_{p}^{2}\\
&  = \tilde{m}_{p}\, t_{p}^{2}+o(1)
\end{align*}
where $o(1)\rightarrow0$ for $p\rightarrow2^{*}.$ Hence the conclusion  follows.
\end{proof}

\section{The barycenter map}	\label{bary}

In this section we introduce the barycenter map that will
allow us to compare the topology of $\Omega$ with the topology of suitable
sublevels of $I_{p}\,;$ precisely sublevels with energy near the minimum level $m_p\,.$

Before to proceed, some other notations are in order.
For $u\in H^{1}(\mathbf{R}^{3})$ with compact support, 
let us denote with the same symbol $u$ its
trivial extension out of supp\,$u$. The barycenter of $u$ (see \cite{BC}) is defined as
\begin{equation*}
\beta(u)=\dfrac{\int_{\mathbf{R}^{3}} x|\nabla u|^{2}\,dx}{\int_{\mathbf{R}
^{3}}|\nabla u|^{2}dx}\,.
\end{equation*}
From now on, we fix  $r>0$ a radius sufficiently small such that $B_{r}\subset\Omega$ and the sets
$$
\Omega^{+}_{r}=\{x\in\mathbf{R}^{3}:d(x,\Omega)\le r\}\,
$$
$$
\Omega^{-}_{r}=\{x\in\Omega:d(x,\partial\Omega)\ge r\}\,
$$
are homotopically equivalent to $\Omega $. In particular we denote by
\begin{equation}\label{h}
h:\Omega^{+}_r \rightarrow\Omega^-_r
\end{equation}
the homotopic equivalence map such that $h|_{\Omega^-_r}$ is the identity. 

Let us  introduce the space $D^{1,2}(\mathbf{R}^{3})=\left\{u\in L^{2^{*}}%
(\mathbf{R}^{3}): \nabla u\in L^{2}\right\}$ which can also be characterized as the
closure of $C^{\infty}_{0}(\mathbf{R}^{3})$ with respect to the (squared)
norm
$$\|u\|^{2}_{D^{1,2}(\mathbf{R}^{3})}=\int_{\mathbf{R}^{3}} |\nabla u|^{2}\,dx.$$
A function in  $H^{1}_{0}(\Omega)$ can be thought as an element of $D^{1,2}(\mathbf{R}^{3})$.

\medskip

The following ``global compactness'' result is taken from Struwe (see Theorem
3.1 of \cite{Stw}) and will be useful to study the behavior of the PS sequences
for the limit functional $I_{*}(u)=\frac{1}{2}\|u\|^{2}-\frac{1}{2^{*}}|u|_{2^{*}}^{2^{*}}$. 

\begin{theorem}\label{St}
Let $\{v_{n}\}$ be a PS sequence for $I_{*}$ in $H^{1}_{0}(\Omega).$
Then there exist a number $k\in\mathbf{N}_{0}$, sequences of points
$\{x_{n}^{j}\} \subset\Omega$ and sequences of radii $\{ R^{j}_{n}\}$ ($1\le
j\le k$) with $R_{n}^{j}\rightarrow+\infty$ for $n\rightarrow+\infty$, there exist a positive solution
$v\in H^{1}_{0}(\Omega)$ of (\ref{star}) and non trivial solutions
$v^{j}\in D^{1,2}(\mathbf{R}^{3} )$ ($1\le
j\le k$) of
\begin{equation}\label{limit}	
	-\Delta u=|u|^{2^{*}-2}\ \ \ \text{ in } \mathbf{R}^{3}\,,
\end{equation}
such that, a (relabeled) subsequence $\{v_{n}\}$ satisfies
\begin{eqnarray*}
&v_{n}-v-\sum_{j=1}^{k} v_{R_{n}}^{j}(\cdot-x_{n}^{j}) \rightarrow0 \ \ \text{
in }\ \ {D^{1,2}(\mathbf{R}^{3})}\,,& \\
&I_{*}(v_{n})\rightarrow I_{*}(v)+\sum_{j=1}^{k}\hat{I}(v^{j})&
\end{eqnarray*}
where $\hat{I}: H^{1}_{0}(\mathbf{R}^{3})\rightarrow\mathbf{R}$ is given by
$$
\hat{I}(u)=\frac{1}{2}\int_{\mathbf{R}^{3}}|\nabla u|^{2} dx-\frac{1}{2^{*}%
}\int_{\mathbf{R}^{3}}|u|^{2^{*}}\,dx\,.
$$
\end{theorem}
Basically the theorem states that if the PS condition fails, it is due to the solutions of \eqref{limit}.
For what concerns $\hat{I}$, it is known that it achieves its minimum on
functions of type
\begin{equation}\label{family}	
	U_{R}(x-a)=\frac{(3 R^{2})^{1/4}}{(R^{2}+|x-a|^{2})^{1/2}}\ \ \ \ R>0\,,a\in \mathbf R^3
\end{equation}
and the  minimum value is exactly $\hat{I}(U_{R}(\cdot-a))=\frac{1}{3}\int_{\mathbf{R}^{3}}|\nabla
U|^{2}dx=m_{*},$ namely the infimum of $I^*.$
On the other hand, the value of $\hat{I}$ on solutions of \eqref{limit} which do not
belong to the family \eqref{family} is greater than $2m_*$.
As a consequence, if the sequence $\{v_n\}$ of Theorem \ref{St} is a PS sequence for $I_*$ at level $m_*$,
we deduce  $I_*(v)=0, k=1$ and $v^{1}=U$. 
Furthermore, since $v$  is a solution of \eqref{star}
and $I_*$ is positive on the solutions,  necessarily $v=0$
and so Theorem \ref{St} gives 
$$v_n-U_{R_n}(\cdot-x_n)\rightarrow 0\ \ \ \mbox{in} \ \ D^{1,2}(\mathbf R^3).$$

Thanks to the previous theorem we can prove that, roughly speaking, if $p$ is near the critical
exponent $2^*$, the functions with barycenter outside $\Omega$ have an energy
away from the ground state level $m_p$.

\begin{proposition}\label{baricentri}
There exists $\varepsilon>0$ such that if $p\in(2^{*}-\varepsilon,2^{*}),$ it
follows
$$u\in\mathcal{N}_{p} \ \mbox{ and } \ I_{p}(u)<m_{p}+\varepsilon\,
\Longrightarrow\,\beta(u)\in\Omega^{+}_{r}.$$
\end{proposition}
\begin{proof}
We argue by contradiction. Assume that there exist sequences $\varepsilon
_{n}\rightarrow0,p_{n}\rightarrow2^{*}$ and $u_{n}\in\mathcal{N}_{p_{n}}$ such
that
\begin{equation}
\label{contradiction}I_{p_{n}}(u_{n})\le m_{p_{n}}+\varepsilon_{n}
\  \mbox{ and } \ \ \beta(u_{n})\notin\Omega^{+}_{r}.
\end{equation}
Then, by Proposition \ref{limitmp}
\begin{equation}\label{converg}	
	I_{p_{n}}(u_{n})\rightarrow m_{*}
\end{equation}
and $\{u_n\}$ is bounded in $\H$. Let $t_{n}>0$ such that $t_{n} u_{n}\in\mathcal{N}_{*}$. By Lemma
\ref{limsuptp} we may assume (up to subsequence) that $t_{n}\rightarrow1$ and we evaluate
\begin{align*}
I_{p_{n}}(u_{n})-I_{*}(t_{n} u_{n}) & =\left( \frac{1}{2}-\frac{1}{p_{n}
}\right) \|u_{n}\|^{2} +\lambda\frac{p_{n}-4}{4p_{n}}\int_{\Omega}\phi_{u_{n}}
u_{n}^{2}\,dx-\left( \frac{1}{2}-\frac{1}{2^{*}}\right) t_{n}^{2}
\|u_{n}\|^{2}\\
&  \ge\left( \frac{1}{2}-\frac{1}{p_{n}}\right) \|u_{n}\|^{2}-\left( \frac
{1}{2}-\frac{1}{2^{*}}\right) t_{n}^{2} \|u_{n}\|^{2}\\
& = \left( \frac{1}{2}-\frac{1}{p_{n}}\right) \|u_{n}\|^{2}\left( 1-t_{n}
^{2}\right) -\left( \frac{1}{p_{n}}-\frac{1}{2^{*}}\right) t_{n}^{2}
\|u_{n}\|^{2}\\
&=o(1)
\end{align*}
which gives
\begin{align*}
m_{*}\le I_{*}(t_{n} u_{n}) & \le I_{p_{n}}(u_{n})+o(1).
\end{align*}
By \eqref{converg}, $I_{*}(t_{n} u_{n})\rightarrow m_{*}$ for $n\rightarrow+\infty.$
The Ekeland's variational principle implies that there exist $\{v_{n}\}\subset
\mathcal{N}_{*} $ and $\{\mu_{n}\}\subset\mathbf{R}$ such that
\begin{align*}
& \|t_{n} u_{n} -v_{n}\|\rightarrow0\\
& I_{*}(v_{n})=\frac{1}{3}\|v_{n}\|^{2}\rightarrow m_{*}\\
& I_{*}^{\prime}(v_{n})-\mu_{n} G_{*}^{\prime}(v_{n})\rightarrow0\nonumber
\end{align*}
and Lemma \ref{lemmaPS} (in the case $\lambda=0$) ensures that $\{v_{n}\}$ is a PS sequence for the free functional
$I_{*}$ at level $m_*$.
By the remarks after Theorem \ref{St},
$$
v_{n} - U_{R_{n}}(\cdot-x_{n})\rightarrow0\ \ \ \text{ in } D^{1,2}
(\mathbf{R}^{3})
$$
where $\{x_{n}\}\subset\Omega, R_{n}\rightarrow+\infty$ and we can write
$$v_{n}=U_{R_{n}}(\cdot-x_{n})+w_{n}$$
with a remainder $w_{n}$ such that $\|w_{n}\|_{D^{1,2}(\mathbf{R}^{3}
)}\rightarrow0 $\,.
It is clear that $t_{n} u_{n}=v_{n}+ t_{n} u_{n}-v_{n}$\,; so, renaming the
remainder again $w_{n}$, we have
$$t_{n} u_{n}=U_{R_{n}}(\cdot-x_{n})+w_{n}.$$
Now writing $x\in\mathbf{R}^{3}$ as $x=(x^{1},x^{2},x^{3})$, the $i-$th
coordinate of the barycenter of $u_{n}$ satisfies
\begin{multline}\label{bari} 
\beta(u_{n})^{i} \| t_{n} u_{n}\|^{2}_{D^{1,2}(\mathbf{R}^{3})}
= \int_{\mathbf{R}^{3}}x^{i} |\nabla U_{R_{n}}(x-x_{n})|^{2}\,dx\\ 
            +\int_{\mathbf{R}^{3}} x^{i}|\nabla w_{n}(x)|^{2}\,dx 
+2\int_{\mathbf{R}^{3}}x^{i} \nabla U_{R_{n}}(x-x_{n})\nabla w_{n}(x)\,dx. 
\end{multline}
The aim is to localize the sequence of barycenters, so we pass to the limit in the above expression  
evaluating $\| t_{n} u_{n}\|^{2}_{D^{1,2}(\mathbf{R}^{3})}$ 
and the right hand side.

First,
\begin{equation}
\label{infm1}\| t_{n} u_{n}\|^{2}_{D^{1,2}(\mathbf{R}^{3})}=\|U\|^{2}%
_{D^{1,2}(\mathbf{R}^{3})}+o(1)\,
\end{equation}
and simple computations show that
\begin{multline}
\label{infm2}\int_{\mathbf{R}^{3}} x^{i} |\nabla U_{R_{n}}(x-x_{n})|^{2}\,dx
=\\
\frac{1}{R_{n}} \int_{\mathbf{R}^{3}} y^{i} \left| \nabla U(y)\right| ^{2}\,dy
+x_{n}^{i} \int_{\mathbf{R}^{3}}|\nabla U(y)|^{2}\,dy\,.
\end{multline}
Moreover, since $v_n$ are supported in $\Omega$,  there holds
\begin{equation*}
U_{R_{n}}(\cdot-x_{n})= - w_n \ \ \text{ on } \mathbf R^3\setminus\Omega
\end{equation*}
and we  evaluate
\begin{align*}
\int_{\mathbf{R}^{3}} x^{i}|\nabla w_{n}(x)|^{2}\,dx & =\int_{\Omega}x^{i}|\nabla
w_{n}(x)|^{2}\,dx+A_n
\end{align*}
where
\begin{align*}
A_{n} & =\int_{\mathbf{R}^{3}\setminus\Omega}x^{i}|\nabla w_{n}(x)|^{2}\,dx\\
&=
\int_{\mathbf{R}^{3}\setminus\Omega} x^{i} R_{n} |\nabla
U(R_{n}(x-x_{n}))|^{2}\,dx\\
& =\int_{\mathbf{R}^{3}\setminus R_{n} (\Omega-x_{n})}(\frac{y^{i}}{R_{n}%
}+x_{n}^{i})|\nabla U(y)|^{2}\,dy\\
& =\frac{1}{R_{n}}\int_{\mathbf{R}^{3}\setminus R_{n} (\Omega-x_{n})}y^{i}
|\nabla U(y)|^{2}\,dy +x_{n}^{i}\int_{\mathbf{R}^{3}\setminus R_{n}(
\Omega-x_{n})}|\nabla U(y)|^{2}\,dy=o(1).
\end{align*}
As a consequence,
\begin{equation}
\label{infm3}\int_{\mathbf{R}^{3}} x^{i}|\nabla w_{n}(x)|^{2}\,dx=\int_{\Omega}
x^{i}|\nabla w_{n}(x)|^{2}\,dx+o(1)=o(1)\,.
\end{equation}
The last term in \eqref{bari} is  estimated as
\begin{multline*}
\int_{\mathbf{R}^{3}} x^{i} \nabla U_{R_{n}}(x-x_{n})\nabla w_{n}(x)\,dx  
= \int_{\Omega}x^{i} \nabla U_{R_{n}}(x-x_{n})\nabla w_{n} (x)\,dx -A_{n}\\
 \le c \left(\int_{\Omega}|\nabla U_{R_n}(x-x_n)|^2\,dx\right)^{1/2} \left(\int_{\Omega}|\nabla w_n|^2\,dx\right)
^{1/2}-A_n 
\end{multline*}
with $A_n$ defined as before and then,
\begin{equation}\label{infm4}
\int_{\mathbf{R}^{3}} x^{i} \nabla U_{R_{n}}(x-x_{n})\nabla w_{n}(x)\,dx =o(1)\,.
\end{equation}
Putting together \eqref{infm1},\eqref{infm2},\eqref{infm3} and \eqref{infm4} by \eqref{bari}
we deduce
\begin{equation}
\label{b}\beta(u_{n})^{i}= \frac{x_{n}^{i} \int_{\mathbf{R}^{3}}|\nabla
U(y)|^{2}\,dy+o(1)}{\|U\|^{2}_{D^{1,2}(\mathbf{R}^3)}+o(1)}\,.
\end{equation}
Since $\{x_{n}\}\subset\Omega$, \eqref{b} implies that definitively $\beta(u_{n})\in\bar{\Omega
}$ which is in contrast with \eqref{contradiction} and proves the proposition.
\end{proof}

\section{Proof of Theorem \ref{Main}}
Here we complete the  proof of our theorem but first we need a slight modi\-- fication to the
previous notations.
We  add a subscript $r$ ($r>0$ and small as before) to denote the same
quantities defined in the previous sections when the domain
$\Omega$ is replaced by $B_{r}$; namely integrals are taken on $B_{r}$ and norms are taken
for functional spaces defined on $B_{r}.$
Hence
$$
\mathcal{N}_{p,r}=\left\{ u\in H^{1}_{0}(B_{r}):\| u\|^{2}_{H^{1}_{0}(B_{r}
)}+\lambda\int_{B_{r}}\phi_{u} u^{2}\,dx=|u|^{p}_{L^{p}(B_{r})}\right\}
$$
and, for $u\in\mathcal{N}_{p,r}$
$$
I_{p,r}(u)=\frac{p-2}{2p}\|u\|^{2}_{H^{1}_{0}(B_{r})}+\lambda\frac{p-4}%
{4p}\int_{B_{r}} \phi_{u} u^{2} dx\,,
$$
$$ m_{p,r}=\min_{\mathcal{N}_{p,r}} I_{p,r}=I_{p,r}(u_{p,r})\,.$$
Moreover let
$$
I_{p}^{m_{p,r}}=\left\{ u\in\mathcal{N}_{p}: I_{p}(u)\le m_{p,r}\right\}
$$
which is non vacuous since $m_{p}<m_{p,r}$.

Define also, for $p\in(4,2^*)$ the map $\Psi_{p,r}:\Omega_{r}^{-}\rightarrow {\mathcal N}_{p}$ such that
\[
\Psi_{p,r}(y)(x)=
\left\{
\begin{array}
[c]{ccl}
u_{p,r}(\left\vert x-y\right\vert ) & \mbox{if} & x\in B_{r}(y)\\
0 & \mbox{if} & x\in\Omega\setminus B_{r}(y)
\end{array}
\right.
\]
and note that we have
\begin{equation*}
	\beta(\Psi_{p,r}(y))=y\ \ \ \mbox{and}\ \ \ \Psi_{r,p}(y)\in I_p^{m_{p,r}}\,.
\end{equation*}
Moreover, 
since $m_p+k_p=m_{p,r}$ where $k_p>0$ and tends to zero if $p\rightarrow2^*$ (see Proposition \ref{limitmp}), in
correspondence of $\varepsilon>0$ provided by Proposition \ref{baricentri}, there exists a $\bar{p}\in [4,2^*)$ such
that for every $p\in[\bar{p},2^*)$ it results $k_p<\varepsilon$; so if $u\in I_p^{m_{p,r}}$ we have
$$I_p(u)\le m_{p,r}<m_p+\varepsilon,$$
at least for $p$ near $2^*.$
Hence the following maps are well-defined:
$$\Omega^-_r\stackrel{\Psi_{p,r}}{\longrightarrow}I_p^{m_{p,r}}\stackrel{h\circ\beta}{\longrightarrow}\Omega_r^-\,$$
where $h$ is given by \eqref{h}.
Since the composite map $ h\circ\beta\circ\Psi_{p,r}$ is the identity of $\Omega^-_r$, 
by a property of the category, 
the sublevel $I_p^{m_{p,r}}$ ``dominates'' the set $\Omega_r^-$ in the sense that
$$\mbox{cat}_{I_p^{m_{p,r}}}(I_p^{m_{p,r}})\ge \mbox{cat}_{\Omega^-_r} (\Omega^-_r)$$
(see e.g. \cite{J}) and our choice of $r$ gives 
$\mbox{cat}_{\Omega^-_r} (\Omega^-_r)=\mbox{cat}_{\bar{\Omega}}(\bar{\Omega})$.
In conclusion, we have found a sublevel of $I_p$ on $\cal N_p$ with category greater than 
$\mbox{cat}_{\bar{\Omega}}(\bar{\Omega})$. 
Since, as we have already said, the PS condition is verified on $\cal N_p$\,, applying the Lusternik-Schnirel\-mann 
theory we get
the existence of at least $\mbox{cat}_{\bar{\Omega}}(\bar{\Omega})$ critical points for $I_p$
on the manifold $\cal N_p$ which give rise to solutions of \eqref{SPS}.

The existence of another solution is obtained with the same arguments of \cite{BC2}. Since by hypothesis
$\Omega$ is not contractible in itself, by the choice of $r$ it results $\mbox{cat}_{\Omega_{r}^{+}}\,(\Omega_{r}^{-})>1$, namely
$\Omega_{r}^{-}$ is not contractible in $\Omega_{r}^{+}$. We claim now that 
 the set ${\Psi_{p,r}(\Omega_{r}^{-})}$ can not be contractible in $I_{p}^{m_{p,r}}$. 
Indeed, assume by contradiction that $\mbox{cat}_{I_{p}^{m_{p,r}}}\,(\Psi_{p,r}(\Omega_{r}^{-}))=1$:
this means that there exists a map $\mathcal H \in C([0,1]\times{\Psi_{p,r}(\Omega_{r}^{-})}; I_{p}^{m_{p,r}})$ such that
$${\mathcal H}(0,u)=u \ \ \forall u\in{\Psi_{p,r}(\Omega_{r}^{-})}\ \  \text{and}$$
$$\exists\, w\in I_{p}^{m_{p,r}}: {\mathcal H}(1,u)=w \ \ \forall u\in {\Psi_{p,r}(\Omega_{r}^{-})} .$$
Then $F=\Psi_{p,r}^{-1}({\Psi_{p,r}(\Omega_{r}^{-})})$ is closed, contains $\Omega_{r}^{-}$ and 
is contractible in $\Omega_{r}^{+}$ as we can see by defining the map 
\begin{equation*}
{\mathcal G}(t,x)=
\begin{cases}
   { \beta(\Psi_{r,p}(x))} & \text{if  $0\le t\le 1/2$}, \\ 
    \beta ({\mathcal H}(2t-1, \Psi_{p,r}(x))) &\text{if  $1/2\le t\le1$}.
\end{cases}
\end{equation*}
Then also $\Omega_{r}^{-}$ would be contractible in $\Omega_{r}^{+}$ giving a contraddiction.

On the other hand we can choose a function $z\in {\mathcal N}_{p}\setminus{\Psi_{p,r}(\Omega_{r}^{-})}$ so that the cone 
$$\mathcal C=\left\lbrace \theta z+(1-\theta) u : u\in {\Psi_{p,r}(\Omega_{r}^{-})}, \theta\in [0,1]\right\rbrace $$
is compact and contractible in $H^{1}_{0}(\Omega)$ and  $0\notin{\mathcal C}$. 
Denoting with $t_{u}$ the unique positive number provided by Lemma \ref{lemmanehari}, it follows that if we set
$$\hat{\mathcal C}=\{t_{u}u : u\in \mathcal C\}, \ \ \ M_{p}=\max_{\hat{\mathcal C}} I_{p}$$
then $\hat{\mathcal C}$ is contractible in $I_{p}^{M_{p}}$ and $M_p>m_{p,r}.$ As a consequence
also ${\Psi_{p,r}(\Omega_{r}^{-})}$ is contractible in $I_{p}^{M_{p}}.$ 

Summing up, the set $\Psi_{p,r}(\Omega_{r}^{-}) $ is contractible in $I_{p}^{M_{p}}$ and not in $I_{p}^{m_{p,r}}$. 
Since the PS condition is satisfied we deduce the existence of another critical point
with critical level between $m_{p,r}$ and $M_{p}$.

It remains to prove that these solutions are positive. Note that we can apply all the previous
machinery replacing the functional \eqref{I} with
\begin{equation*}
I^+_{p}(u)=\frac{1}{2}\int_{\O}(|\nabla u|^{2}+u^{2})\,dx+\frac{\lambda}{4}
\int_{\O}\phi_{u} u^{2}dx-\frac{1}{p}\int_{\O}(u^+)^{p} dx
\end{equation*}
obtaining again at least $\mbox{cat}_{\bar{\Omega}}(\bar{\Omega})$ nontrivial solutions. Finally
the maximum principle ensures that these solutions are positive, hence they solve \eqref{SPS}.

\medskip


\end{document}